\documentclass[12pt]{article}%
\usepackage{amsfonts}
\usepackage{sw20bams}
\usepackage{amsmath}
\usepackage{amssymb}
\usepackage{graphicx}%
\providecommand{\U}[1]{\protect\rule{.1in}{.1in}}
\begin{document}

\title{Small Matrices with Small Inverses: Unimodular Zerofree Cases}
\author{Steven Finch}
\date{June 12, 2026}
\maketitle

\begin{abstract}
We consider unimodular matrices $M$ such that neither $M$ nor $M^{-1}$ contain
zero entries. Matrices typically exhibit a trade-off: small $M$ imply large
$M^{-1}$. We investigate rare cases where both remain small, classify these
matrices up to symmetry, and discuss aspects of this balanced setting.

\end{abstract}

\footnotetext{Copyright \copyright \ 2026 by Steven R. Finch. All rights
reserved.}

An $n\times n$ integer matrix $M$ with determinant $\pm1$ is called
\textbf{unimodular}. \ Let $\alpha=\left\Vert M\right\Vert $ and
$\beta=\left\Vert M^{-1}\right\Vert $, the maximum absolute entry of $M$ and
$M^{-1}$, respectively. \ Examples:\ $\alpha=\beta=1$ for the identity matrix
$I$ or its negative $-I$, and $\alpha=\beta=2$ for%
\[%
\begin{array}
[c]{ccc}%
M_{3\times3}=\left(
\begin{array}
[c]{ccc}%
1 & 0 & 0\\
0 & 1 & 1\\
0 & 1 & 2
\end{array}
\right)  , &  & M_{5\times5}=\left(
\begin{array}
[c]{ccccc}%
1 & 0 & 0 & 0 & 0\\
0 & 1 & 1 & 0 & 0\\
0 & 1 & 2 & 0 & 0\\
0 & 0 & 0 & 1 & 1\\
0 & 0 & 0 & 1 & 2
\end{array}
\right)  .
\end{array}
\]
The matrix $M$ is called \textbf{zerofree} if none of the entries in $M$ and
none of the entries in $M^{-1}$ are zero. \ (Involving both matrices in this
definition is somewhat unorthodox.) \ Our interest is in unimodular zerofree
matrices for which both $\alpha$ and $\beta$ are small.\smallskip

\textit{Proposition 0}. \ If $n>1$, then $\alpha>1$ and $\beta>1$.\smallskip

\textit{Proof}. \ Suppose that $M$ satisfies $\alpha=1$. \ Clearly $M\equiv
J\operatorname{mod}2$ where $J$ is the matrix of ones, since each $m_{ij}%
\in\{-1,1\}$. \ For $n>1$, $J$ has rank $1$ over $\mathbb{F}_{2}$, thus
$\det(M)\equiv\det(J)\equiv0\operatorname{mod}2$. \ It follows that $\det(M)$
is even, which contradicts $\det(M)=\pm1$. \ Replace $M$ by $M^{-1}$ to
similarly exclude $\beta=1$. This simple proof, although well-known,
apparently does not appear in the published literature. \ In fact, $\det(M)$
is divisible by $2^{n-1}$ via a more complicated argument \cite{Cohn-mtrx2,
MWS-mtrx2, DP-mtrx2, Zivk-matrx2}.\smallskip

We wish to classify unimodular zerofree matrices (for fixed $n$, $\alpha$,
$\beta$) up to the double action of the signed-permutation group, so that two
matrices $M$, $N$ are equivalent when $N=P\,M\,Q$ for some matrices $P$ \& $Q$
that independently rearrange rows, rearrange columns, and flip the signs of
rows \& columns. \ This action partitions the space into \textbf{orbits}, each
consisting of all matrices obtainable from one another by such
transformations. To assign each orbit a unique \textbf{canonical
representative}, we compare matrices by their row-major flattenings, sorted
lexicographically using the structural integer ordering%
\[
1<2<3<4<5<\cdots<-1<-2<-3<-4<-5<\cdots
\]
and we select the matrix whose flattening is minimal in this ordering. This
canonical form provides an unambiguous label for each equivalence class and
allows orbit membership to be detected by a single comparison.

Henceforth, when we refer to a \textquotedblleft matrix\textquotedblright, we
mean the canonical representative of its equivalence class, and we freely
identify each matrix with its orbit under the signed-permutation action.
\ Also, when we describe a matrix as \textquotedblleft
positive\textquotedblright, we mean that its canonical representative has all
positive entries; other members of the orbit might not share this property.

\section{Cases When $n=2$}

For $\alpha=\beta=2$, there is a unique matrix%
\[
\left(
\begin{array}
[c]{cc}%
1 & 1\\
1 & 2
\end{array}
\right)  .
\]
For $\alpha=\beta=3$, three matrices exist:%
\[%
\begin{array}
[c]{ccccc}%
\left(
\begin{array}
[c]{cc}%
1 & 2\\
2 & 3
\end{array}
\right)  , &  & \left(
\begin{array}
[c]{cc}%
1 & 2\\
1 & 3
\end{array}
\right)  , &  & \left(
\begin{array}
[c]{cc}%
1 & 1\\
2 & 3
\end{array}
\right)
\end{array}
\]
and likewise for $\alpha=\beta=4$:%
\[%
\begin{array}
[c]{ccccc}%
\left(
\begin{array}
[c]{cc}%
2 & 3\\
3 & 4
\end{array}
\right)  , &  & \left(
\begin{array}
[c]{cc}%
1 & 3\\
1 & 4
\end{array}
\right)  , &  & \left(
\begin{array}
[c]{cc}%
1 & 1\\
3 & 4
\end{array}
\right)  .
\end{array}
\]
For $\alpha=\beta=5$, seven matrices exist:%
\[%
\begin{array}
[c]{ccccccc}%
\left(
\begin{array}
[c]{cc}%
3 & 4\\
4 & 5
\end{array}
\right)  , & \left(
\begin{array}
[c]{cc}%
2 & 3\\
3 & 5
\end{array}
\right)  , & \left(
\begin{array}
[c]{cc}%
1 & 4\\
1 & 5
\end{array}
\right)  , & \left(
\begin{array}
[c]{cc}%
1 & 3\\
2 & 5
\end{array}
\right)  , & \left(
\begin{array}
[c]{cc}%
1 & 2\\
3 & 5
\end{array}
\right)  , & \left(
\begin{array}
[c]{cc}%
1 & 2\\
2 & 5
\end{array}
\right)  , & \left(
\begin{array}
[c]{cc}%
1 & 1\\
4 & 5
\end{array}
\right)
\end{array}
\]
but only three for $\alpha=\beta=6$:%
\[%
\begin{array}
[c]{ccccc}%
\left(
\begin{array}
[c]{cc}%
4 & 5\\
5 & 6
\end{array}
\right)  , &  & \left(
\begin{array}
[c]{cc}%
1 & 5\\
1 & 6
\end{array}
\right)  , &  & \left(
\begin{array}
[c]{cc}%
1 & 1\\
5 & 6
\end{array}
\right)  .
\end{array}
\]
Counts of matrices for $2\leq\alpha=\beta<\infty$ constitute an interesting
sequence%
\[
1,3,3,7,3,11,7,11,7,19,7,23,11,15,15,31,11,35,15,23,19,43,15,39,23,35,23,55,15,\ldots
\]
worthy of further study. \ [See the Addendum.]

\section{Cases When $n=3$}

For $\alpha=\beta=3$, there is a unique matrix%
\[
\left(
\begin{array}
[c]{ccc}%
1 & 2 & 2\\
2 & 1 & 2\\
2 & 2 & 3
\end{array}
\right)  .
\]
For $\alpha=3<\beta=4$ \& $\alpha=4>\beta=3$, there are unique matrices%
\[%
\begin{array}
[c]{ccc}%
\left(
\begin{array}
[c]{ccc}%
1 & 1 & 1\\
1 & 2 & -1\\
2 & 3 & -1
\end{array}
\right)  , &  & \left(
\begin{array}
[c]{ccc}%
1 & 1 & 1\\
1 & 2 & 3\\
1 & 3 & 4
\end{array}
\right)
\end{array}
\]
respectively, and likewise for $\alpha=2<\beta=5$ \& $\alpha=5>\beta=2$:%
\[%
\begin{array}
[c]{ccc}%
\left(
\begin{array}
[c]{ccc}%
1 & 1 & 2\\
1 & -2 & -2\\
2 & -2 & -1
\end{array}
\right)  , &  & \left(
\begin{array}
[c]{ccc}%
2 & 2 & 3\\
2 & 3 & 4\\
3 & 4 & 5
\end{array}
\right)  .
\end{array}
\]
The left-hand matrices are canonical representatives of inverses of the
right-hand matrices. \ A careless tendency to treat $\alpha\leq\beta$ as valid
\textquotedblleft without loss of generality\textquotedblright\ is
demonstrably wrong. \ While the following can be verified computationally, a
conceptual proof is not known.\smallskip

\textit{Conjecture 1}. \ \ No $3\times3$ matrices exist that are tied to
$\alpha=2$ and $2\leq\beta\leq4$.\smallskip

\noindent For $\alpha=3<\beta=5$, six matrices exist:%
\[%
\begin{array}
[c]{ccccc}%
\left(
\begin{array}
[c]{ccc}%
1 & 2 & 2\\
3 & 1 & 2\\
3 & 2 & 3
\end{array}
\right)  , &  & \left(
\begin{array}
[c]{ccc}%
1 & 2 & 2\\
2 & 2 & 3\\
3 & 1 & 3
\end{array}
\right)  , &  & \left(
\begin{array}
[c]{ccc}%
1 & 1 & 2\\
1 & 2 & 3\\
1 & -2 & -2
\end{array}
\right)  ,
\end{array}
\]%
\[%
\begin{array}
[c]{ccccc}%
\left(
\begin{array}
[c]{ccc}%
1 & 1 & 1\\
1 & 2 & 3\\
1 & -1 & -2
\end{array}
\right)  , &  & \left(
\begin{array}
[c]{ccc}%
1 & 1 & 1\\
1 & 2 & -1\\
1 & 3 & -2
\end{array}
\right)  , &  & \left(
\begin{array}
[c]{ccc}%
1 & 1 & 1\\
1 & 2 & -2\\
2 & 3 & -2
\end{array}
\right)
\end{array}
\]
and likewise for $\alpha=\beta=4$:%
\[%
\begin{array}
[c]{ccccc}%
\left(
\begin{array}
[c]{ccc}%
1 & 3 & 3\\
2 & 2 & 3\\
2 & 3 & 4
\end{array}
\right)  , &  & \left(
\begin{array}
[c]{ccc}%
1 & 2 & 2\\
3 & 2 & 3\\
3 & 3 & 4
\end{array}
\right)  , &  & \left(
\begin{array}
[c]{ccc}%
1 & 1 & 3\\
1 & 2 & 4\\
2 & 1 & 4
\end{array}
\right)  ,
\end{array}
\]%
\[%
\begin{array}
[c]{ccccc}%
\left(
\begin{array}
[c]{ccc}%
1 & 1 & 2\\
1 & 2 & 3\\
1 & 4 & 4
\end{array}
\right)  , &  & \left(
\begin{array}
[c]{ccc}%
1 & 1 & 2\\
1 & 2 & 1\\
3 & 4 & 4
\end{array}
\right)  , &  & \left(
\begin{array}
[c]{ccc}%
1 & 1 & 1\\
1 & 2 & 4\\
2 & 3 & 4
\end{array}
\right)  .
\end{array}
\]
For $\alpha=3<\beta=6$, seven matrices exist:%
\[%
\begin{array}
[c]{ccccccc}%
\left(
\begin{array}
[c]{ccc}%
1 & 1 & 2\\
1 & 3 & 1\\
2 & 3 & 3
\end{array}
\right)  , &  & \left(
\begin{array}
[c]{ccc}%
1 & 1 & 2\\
1 & 2 & 3\\
3 & 1 & 3
\end{array}
\right)  , &  & \left(
\begin{array}
[c]{ccc}%
1 & 1 & 2\\
1 & -3 & -3\\
2 & -3 & -2
\end{array}
\right)  , &  & \left(
\begin{array}
[c]{ccc}%
1 & 1 & 2\\
1 & 2 & 3\\
1 & -3 & -3
\end{array}
\right)  ,
\end{array}
\]%
\[%
\begin{array}
[c]{ccccc}%
\left(
\begin{array}
[c]{ccc}%
1 & 2 & 2\\
1 & -2 & -3\\
2 & -1 & -2
\end{array}
\right)  , &  & \left(
\begin{array}
[c]{ccc}%
1 & 1 & 2\\
2 & -2 & 1\\
3 & -2 & 2
\end{array}
\right)  , &  & \left(
\begin{array}
[c]{ccc}%
1 & 1 & 1\\
1 & 2 & -3\\
2 & 3 & -3
\end{array}
\right)
\end{array}
\]
but only four for $\alpha=4<\beta=5$:%
\[%
\begin{array}
[c]{ccccccc}%
\left(
\begin{array}
[c]{ccc}%
2 & 2 & 3\\
2 & 3 & 4\\
3 & 2 & 4
\end{array}
\right)  , &  & \left(
\begin{array}
[c]{ccc}%
2 & 2 & 3\\
2 & 3 & 2\\
3 & 4 & 4
\end{array}
\right)  , &  & \left(
\begin{array}
[c]{ccc}%
1 & 2 & 3\\
2 & 3 & 4\\
3 & 4 & 4
\end{array}
\right)  , &  & \left(
\begin{array}
[c]{ccc}%
1 & 1 & 1\\
1 & -2 & -3\\
1 & -3 & -4
\end{array}
\right)  .
\end{array}
\]
We had thought that a pattern might be inferred from $\alpha+\beta=7$ (one
matrix) and $\alpha+\beta=8$ (six matrices), but this behavior failed to carry
over to $\alpha+\beta=9$. \ 

\section{Cases When $n=4$}

For $\alpha=\beta=2$, three matrices exist:%
\[%
\begin{array}
[c]{ccccc}%
\left(
\begin{array}
[c]{cccc}%
1 & 1 & 1 & 2\\
1 & 1 & 2 & 1\\
1 & 2 & 2 & 2\\
2 & 1 & 2 & 2
\end{array}
\right)  , &  & \left(
\begin{array}
[c]{cccc}%
1 & 1 & 1 & 2\\
1 & 2 & 2 & 2\\
1 & -1 & -2 & 1\\
2 & -1 & -2 & 2
\end{array}
\right)  , &  & \left(
\begin{array}
[c]{cccc}%
1 & 1 & 1 & 1\\
1 & 1 & 2 & 2\\
1 & 2 & 1 & 2\\
1 & 2 & -1 & 1
\end{array}
\right)  .
\end{array}
\]
It is fascinating that $\alpha=\beta=2$ was impossible in the preceding
section. \ While the following can be verified computationally, a conceptual
proof is not known.\smallskip

\textit{Conjecture 2}. \ \ No $4\times4$ matrices exist that are tied to
$\alpha=2$ and $\beta=3$.\smallskip

\noindent For $\alpha=2<\beta=4$, we have a unique matrix%
\[
\left(
\begin{array}
[c]{cccc}%
1 & 1 & 1 & 1\\
1 & 1 & 2 & 2\\
1 & -1 & 1 & 2\\
1 & -2 & -1 & 1
\end{array}
\right)  .
\]
For $\alpha=2<\beta=5$, six matrices exist:%
\[%
\begin{array}
[c]{ccccc}%
\left(
\begin{array}
[c]{cccc}%
1 & 1 & 1 & 2\\
1 & 1 & 2 & -1\\
1 & 2 & 2 & 1\\
2 & 1 & 2 & 1
\end{array}
\right)  , &  & \left(
\begin{array}
[c]{cccc}%
1 & 1 & 1 & 2\\
1 & 2 & 2 & 2\\
1 & -1 & -2 & -2\\
2 & -1 & -2 & -1
\end{array}
\right)  , &  & \left(
\begin{array}
[c]{cccc}%
1 & 1 & 1 & 2\\
1 & 1 & 2 & 1\\
1 & 2 & 2 & 2\\
2 & -1 & 1 & 1
\end{array}
\right)  ,
\end{array}
\]%
\[%
\begin{array}
[c]{ccccc}%
\left(
\begin{array}
[c]{cccc}%
1 & 1 & 1 & 1\\
1 & 1 & 2 & 2\\
1 & 2 & 1 & 2\\
1 & -1 & -1 & -2
\end{array}
\right)  , &  & \left(
\begin{array}
[c]{cccc}%
1 & 1 & 1 & 2\\
1 & 2 & 2 & 2\\
1 & -1 & -2 & 1\\
2 & 2 & 1 & 2
\end{array}
\right)  , &  & \left(
\begin{array}
[c]{cccc}%
1 & 1 & 1 & 1\\
1 & 1 & 2 & -1\\
1 & 2 & 1 & -1\\
1 & 2 & 2 & -2
\end{array}
\right)
\end{array}
\]
and for $\alpha=2<\beta=6$, six matrices again exist:%
\[%
\begin{array}
[c]{ccccc}%
\left(
\begin{array}
[c]{cccc}%
1 & 1 & 2 & 2\\
1 & 2 & 1 & 2\\
1 & 2 & -2 & 1\\
2 & 2 & -1 & 2
\end{array}
\right)  , &  & \left(
\begin{array}
[c]{cccc}%
1 & 1 & 1 & 2\\
1 & 2 & 2 & 2\\
2 & 1 & 2 & 2\\
2 & -2 & 1 & -1
\end{array}
\right)  , &  & \left(
\begin{array}
[c]{cccc}%
1 & 1 & 1 & 2\\
1 & 1 & 2 & 1\\
1 & -1 & 2 & -2\\
1 & -2 & 1 & -2
\end{array}
\right)  ,
\end{array}
\]%
\[%
\begin{array}
[c]{ccccc}%
\left(
\begin{array}
[c]{cccc}%
1 & 1 & 1 & 1\\
1 & 1 & 2 & -1\\
1 & -1 & 1 & -2\\
2 & -1 & 1 & -1
\end{array}
\right)  , &  & \left(
\begin{array}
[c]{cccc}%
1 & 1 & 1 & 1\\
1 & 1 & 2 & 2\\
1 & -2 & 1 & -1\\
2 & -2 & 1 & -2
\end{array}
\right)  , &  & \left(
\begin{array}
[c]{cccc}%
1 & 1 & 1 & 1\\
1 & 1 & 2 & 2\\
1 & -1 & 1 & -2\\
1 & -2 & 2 & -2
\end{array}
\right)  .
\end{array}
\]
Counts of matrices for $\alpha=2<4\leq\beta\leq26$ constitute another
interesting sequence%
\[
1,6,6,12,10,15,10,16,19,16,11,26,14,16,11,12,20,11,12,0,10,0,8.
\]
For $\alpha=\beta=3$, we count $163$ matrices, of which $38$ are positive.

\section{Cases When $n=5$}

For $\alpha=2<\beta=3$, two matrices exist:
\[%
\begin{array}
[c]{ccc}%
\left(
\begin{array}
[c]{ccccc}%
1 & 1 & 1 & 2 & 2\\
1 & 1 & 2 & 1 & 2\\
1 & 2 & 2 & 2 & 2\\
2 & 1 & 2 & 2 & 2\\
2 & 2 & 2 & 2 & 1
\end{array}
\right)  , &  & \left(
\begin{array}
[c]{ccccc}%
1 & 1 & 1 & 1 & 1\\
1 & 1 & 1 & 2 & -1\\
1 & 1 & 2 & 1 & -1\\
1 & 2 & 2 & 2 & -1\\
2 & 1 & 2 & 2 & -1
\end{array}
\right)  .
\end{array}
\]
While the following can be verified computationally, a conceptual proof is not
known.\smallskip

\textit{Conjecture 3}. \ \ No $5\times5$ matrices exist that are tied to
$\alpha=\beta=2$.\smallskip

\noindent For $\alpha=2<\beta=4$, twenty-two matrices exist, including one
positive matrix:
\[%
\begin{array}
[c]{c}%
\left(
\begin{array}
[c]{ccccc}%
1 & 1 & 1 & 2 & 2\\
1 & 1 & 2 & 1 & 2\\
1 & 2 & 1 & 1 & 2\\
2 & 1 & 1 & 2 & 1\\
2 & 2 & 2 & 1 & 2
\end{array}
\right)  ,
\end{array}
\]
matrices with one or two negative entries:
\[%
\begin{array}
[c]{ccc}%
\left(
\begin{array}
[c]{ccccc}%
1 & 1 & 1 & 2 & 2\\
1 & 1 & 2 & 1 & 2\\
1 & 2 & 2 & 2 & 2\\
2 & 1 & 2 & -1 & 1\\
2 & 2 & 2 & 1 & 2
\end{array}
\right)  , &  & \left(
\begin{array}
[c]{ccccc}%
1 & 1 & 1 & 1 & 1\\
1 & 1 & 1 & 2 & 2\\
1 & 1 & 2 & 1 & 2\\
1 & 2 & 2 & -1 & 1\\
1 & -1 & 1 & 2 & 1
\end{array}
\right)
\end{array}
\]
and the remaining matrices given as vectors in row-major order:%
\[%
\begin{tabular}
[c]{|l|}\hline
$(n,\alpha,\beta)=(5,2,4)$ matrices with $3,4,\ldots,9$ negative
entries:\\\hline
{\small 1, 1, 1, 1, 2, 1, 1, 2, 2, 2, 1, -2, 1, 2, 1, 2, -1, 2, 2, 2, 2, -2,
1, 2, 2}\\\hline
{\small 1, 1, 1, 2, 2, 1, 1, 2, 1, 2, 1, 2, 2, 2, 2, 1, -2, 1, -2, -1, 2, 1,
2, 2, 2}\\\hline
{\small 1, 1, 1, 1, 2, 1, 1, 2, 2, 2, 1, 2, 1, 2, 2, 1, 2, -1, 1, 1, 1, 2, -1,
2, -1}\\\hline
{\small 1, 1, 1, 1, 1, 1, 1, 1, 2, 2, 1, 1, 2, 1, -1, 1, 2, 1, 1, -1, 1, 2, 2,
2, -1}\\\hline
{\small 1, 1, 1, 1, 1, 1, 1, 2, 2, 2, 1, 2, 1, 2, 2, 1, 2, -1, 1, -1, 2, 2, 1,
2, -1}\\\hline
{\small 1, 1, 1, 1, 1, 1, 1, 1, 2, 2, 1, 1, 2, 1, 2, 1, 2, 1, 1, 2, 1, 2, -1,
-1, -1}\\\hline
{\small 1, 1, 1, 1, 2, 1, 1, 1, 2, 1, 1, 1, 2, 2, 2, 1, 2, 1, 2, 2, 1, -2, -2,
-2, -2}\\\hline
{\small 1, 1, 1, 1, 2, 1, 1, 1, 2, 1, 1, 1, 2, 2, 2, 1, 2, 1, 2, 2, 1, -1, -2,
-1, -2}\\\hline
{\small 1, 1, 1, 1, 1, 1, 1, 1, 2, -1, 1, 1, 2, 1, -2, 1, 2, 2, 2, -1, 2, 1,
2, 2, -2}\\\hline
{\small 1, 1, 1, 1, 1, 1, 1, 1, 2, -2, 1, 1, 2, 1, -2, 1, 2, 2, 2, -2, 2, 1,
2, 2, -2}\\\hline
{\small 1, 1, 1, 1, 2, 1, 1, 2, -2, -1, 1, 2, 2, -1, 1, 2, 1, 2, -1, 1, 2, 2,
2, -1, 2}\\\hline
{\small 1, 1, 1, 1, 2, 1, 1, 2, 2, 2, 1, 2, -1, -2, 1, 1, -1, 2, 2, 1, 2, 2,
-1, -2, 2}\\\hline
{\small 1, 1, 1, 1, 2, 1, 1, 2, 2, 2, 1, 1, -1, -2, 1, 1, -1, 2, 2, -1, 2, -1,
2, 1, -1}\\\hline
{\small 1, 1, 1, 1, 2, 1, 2, 2, -1, 1, 1, -1, -2, 2, 1, 2, 2, 2, -1, 2, 2, -1,
-2, 1, 1}\\\hline
{\small 1, 1, 1, 1, 1, 1, 1, 1, 2, 2, 1, 2, 2, -1, -2, 1, -1, -2, 1, 1, 2, 2,
1, -1, -2}\\\hline
{\small 1, 1, 1, 1, 2, 1, 1, 2, -2, -1, 1, 2, 2, -1, 1, 2, -1, 2, 1, -1, 2,
-2, 2, 1, -2}\\\hline
{\small 1, 1, 1, 1, 2, 1, 1, 2, 2, 2, 1, -2, 1, -1, -2, 2, -1, 2, 1, 1, 2, -2,
1, -1, -1}\\\hline
{\small 1, 1, 1, 1, 2, 1, 1, 2, 2, 2, 1, -2, -1, -1, -2, 2, -1, 1, 2, -1, 2,
-2, 1, 2, -2}\\\hline
{\small 1, 1, 1, 1, 1, 1, 1, 1, 2, 2, 1, 1, -1, -1, -2, 1, 2, -1, -1, -1, 1,
2, -2, -1, -2}\\\hline
\end{tabular}
\
\]
\noindent For $\alpha=\beta=3$, we count $1352$ matrices, of which $189$ are positive.

\section{Cases When $n=6$}

For $\alpha=\beta=2$, four positive matrices exist:
\[%
\begin{array}
[c]{ccc}%
\left(
\begin{array}
[c]{cccccc}%
1 & 1 & 1 & 1 & 1 & 2\\
1 & 1 & 1 & 1 & 2 & 1\\
1 & 1 & 1 & 2 & 1 & 1\\
1 & 1 & 2 & 1 & 2 & 2\\
1 & 2 & 1 & 2 & 1 & 2\\
2 & 1 & 1 & 2 & 2 & 1
\end{array}
\right)  , &  & \left(
\begin{array}
[c]{cccccc}%
1 & 1 & 1 & 1 & 1 & 2\\
1 & 1 & 1 & 1 & 2 & 1\\
1 & 1 & 1 & 2 & 2 & 2\\
1 & 1 & 2 & 2 & 1 & 2\\
1 & 2 & 2 & 1 & 1 & 2\\
2 & 1 & 2 & 2 & 2 & 2
\end{array}
\right)  ,
\end{array}
\]%
\[%
\begin{array}
[c]{ccc}%
\left(
\begin{array}
[c]{cccccc}%
1 & 1 & 1 & 1 & 1 & 2\\
1 & 1 & 1 & 2 & 2 & 2\\
1 & 1 & 2 & 1 & 2 & 2\\
1 & 2 & 1 & 2 & 2 & 1\\
1 & 2 & 2 & 2 & 2 & 2\\
2 & 2 & 2 & 1 & 2 & 2
\end{array}
\right)  , &  & \left(
\begin{array}
[c]{cccccc}%
1 & 1 & 1 & 2 & 2 & 2\\
1 & 2 & 2 & 1 & 1 & 2\\
1 & 2 & 2 & 2 & 2 & 2\\
2 & 1 & 2 & 1 & 2 & 1\\
2 & 1 & 2 & 2 & 2 & 2\\
2 & 2 & 2 & 1 & 2 & 2
\end{array}
\right)
\end{array}
\]
as well as $199$ matrices containing at least one negative entry.
\ A\ voluminous table of these matrices appears in Addendum II\ of
\cite{Fnch-mtrx2}. \smallskip

\textit{Conjecture 4}. \ \ No $6\times6$ matrices exist that are tied to
$\alpha=\beta=3$.\smallskip

\noindent The cases of $\alpha=2<\beta=3$ and $\alpha=3<4\leq\beta\leq5$ are
similarly complicated. \ 

\section{Cases When $n=7$}

For $\alpha=\beta=2$, two positive matrices exist:%
\[%
\begin{array}
[c]{ccc}%
\left(
\begin{array}
[c]{ccccccc}%
1 & 1 & 1 & 1 & 1 & 2 & 2\\
1 & 1 & 1 & 1 & 2 & 1 & 2\\
1 & 1 & 1 & 1 & 2 & 2 & 1\\
1 & 1 & 1 & 2 & 2 & 2 & 2\\
1 & 1 & 2 & 1 & 2 & 2 & 2\\
1 & 2 & 1 & 1 & 2 & 2 & 2\\
2 & 1 & 1 & 1 & 2 & 2 & 2
\end{array}
\right)  , &  & \left(
\begin{array}
[c]{ccccccc}%
1 & 1 & 1 & 1 & 1 & 1 & 2\\
1 & 1 & 1 & 1 & 1 & 2 & 1\\
1 & 1 & 1 & 1 & 2 & 1 & 1\\
1 & 1 & 1 & 2 & 1 & 1 & 1\\
1 & 2 & 2 & 2 & 2 & 2 & 2\\
2 & 1 & 2 & 2 & 2 & 2 & 2\\
2 & 2 & 1 & 2 & 2 & 2 & 2
\end{array}
\right)  .
\end{array}
\]
Notice that the left-hand matrix is the transpose of the right-hand matrix.
\ There are$\ \ \geq18$ matrices containing at least one negative entry; these
are scantly listed in Addendum III\ of \cite{Fnch-mtrx2}.

For $\alpha=2<\beta=3$, we count $35$ positive matrices:%
\[%
\begin{tabular}
[c]{|l|}\hline
$(n,\alpha,\beta)=(7,2,3)$ matrices with only positive entries [first
part]\\\hline
{\small 1 1 1 1 1 1 2 1 1 1 1 1 2 1 1 1 1 1 2 1 1 1 1 1 2 1 1 1 1 1 2 1 2 2 2
1 2 1 1 2 2 2 2 1 1 1 2 2 2}\\\hline
{\small 1 1 1 1 1 1 2 1 1 1 1 1 2 1 1 1 1 1 2 2 2 1 1 1 2 1 2 2 1 1 2 1 2 1 1
1 2 2 2 1 2 2 2 1 2 2 1 2 2}\\\hline
{\small 1 1 1 1 1 2 2 1 1 1 1 2 1 2 1 1 1 2 1 1 2 1 1 2 1 2 2 2 1 2 1 2 1 2 2
2 1 2 2 2 1 1 2 2 1 2 2 2 2}\\\hline
{\small 1 1 1 1 1 2 2 1 1 1 1 2 1 2 1 1 1 2 1 2 1 1 1 2 2 2 1 2 1 2 1 2 2 1 2
2 1 2 1 1 2 2 2 2 1 2 2 2 2}\\\hline
{\small 1 1 1 1 1 2 2 1 1 1 1 2 1 2 1 1 1 2 1 2 1 1 1 2 1 2 2 2 1 2 1 2 2 2 1
2 1 2 2 2 2 2 2 2 1 2 1 2 1}\\\hline
{\small 1 1 1 1 1 2 2 1 1 1 1 2 1 2 1 1 1 2 1 2 1 1 1 2 1 2 2 2 1 2 1 2 1 2 2
2 1 1 2 2 1 2 2 2 2 2 2 1 2}\\\hline
{\small 1 1 1 1 1 2 2 1 1 1 1 2 1 2 1 1 1 2 1 1 2 1 1 2 1 2 2 2 1 2 1 2 1 2 2
2 1 1 2 2 1 2 2 2 2 2 2 2 1}\\\hline
{\small 1 1 1 1 1 2 2 1 1 1 1 2 1 2 1 1 2 2 1 2 2 1 2 1 2 2 2 1 1 2 2 2 1 2 1
2 2 2 1 2 2 2 2 2 2 2 2 2 1}\\\hline
{\small 1 1 1 1 1 2 2 1 1 1 1 2 1 2 1 1 2 2 1 2 2 1 2 1 2 1 2 2 1 2 2 2 2 2 2
2 1 2 1 2 2 1 2 1 2 2 2 2 2}\\\hline
{\small 1 1 1 1 1 2 2 1 1 1 2 2 2 2 1 1 2 1 2 2 2 1 2 1 2 1 2 2 1 2 2 1 2 1 2
2 2 2 1 1 2 1 2 2 2 1 2 2 2}\\\hline
{\small 1 1 1 1 1 2 2 1 1 1 1 2 1 2 1 1 1 2 2 2 2 1 1 2 2 2 2 1 1 2 2 1 1 2 2
2 1 2 2 2 2 2 2 2 2 2 1 2 2}\\\hline
{\small 1 1 1 1 1 2 2 1 1 1 2 2 2 2 1 1 2 1 2 2 2 1 1 2 2 2 1 2 1 2 1 2 1 2 2
2 1 2 1 2 2 1 2 2 2 2 1 2 2}\\\hline
\end{tabular}
\
\]

\begin{center}%
\[%
\begin{tabular}
[c]{|l|}\hline
$(n,\alpha,\beta)=(7,2,3)$ matrices with only positive entries [second
part]\\\hline
{\small 1 1 1 1 1 2 2 1 1 1 1 2 1 2 1 1 2 2 1 2 2 1 2 1 2 2 1 2 1 2 2 2 1 2 1
2 2 2 1 2 2 2 2 2 2 2 2 1 2}\\\hline
{\small 1 1 1 1 1 2 2 1 1 1 2 2 2 2 1 1 2 1 2 2 2 1 1 2 2 2 1 2 1 2 1 2 1 2 2
2 1 2 2 2 2 2 2 2 2 1 1 2 1}\\\hline
{\small 1 1 1 1 1 2 2 1 1 1 1 2 1 2 1 1 1 2 2 2 2 1 1 2 2 1 2 2 1 2 2 1 1 2 2
2 1 2 2 2 2 2 2 2 2 2 2 2 1}\\\hline
{\small 1 1 1 1 1 2 2 1 1 1 2 2 1 1 1 1 2 1 2 2 2 1 1 2 2 2 1 2 1 2 1 1 2 2 2
2 1 2 2 2 2 2 2 2 2 2 2 1 2}\\\hline
{\small 1 1 1 1 1 2 2 1 1 1 2 2 2 2 1 1 2 2 2 1 2 1 2 2 1 2 1 2 1 2 2 2 2 1 1
2 2 1 1 1 2 2 2 2 2 2 1 2 2}\\\hline
{\small 1 1 1 1 1 2 2 1 1 1 2 2 1 1 1 1 2 1 2 2 2 1 2 2 1 1 2 2 1 2 2 2 1 1 2
2 2 1 2 2 2 2 2 2 2 2 2 1 2}\\\hline
{\small 1 1 1 1 1 2 2 1 1 1 1 2 1 2 1 1 1 2 2 2 2 1 1 2 2 1 2 2 1 2 2 2 1 2 1
2 1 2 2 2 2 2 2 2 2 1 2 2 2}\\\hline
{\small 1 1 1 1 1 2 2 1 1 1 1 2 1 2 1 1 1 2 2 2 2 1 1 2 1 2 2 2 1 2 1 2 2 2 1
2 1 2 2 2 2 2 2 2 1 2 2 2 2}\\\hline
{\small 1 1 1 1 1 2 2 1 1 1 1 2 1 2 1 1 1 2 1 2 1 1 1 2 1 2 2 2 1 2 2 2 2 2 2
2 1 2 2 2 2 2 2 2 2 2 1 2 2}\\\hline
{\small 1 1 1 1 1 2 2 1 1 1 1 2 1 2 1 1 1 2 2 2 2 1 1 2 2 2 2 1 1 2 2 1 2 2 1
2 1 2 2 2 2 2 2 2 2 2 2 2 1}\\\hline
{\small 1 1 1 1 1 2 2 1 1 1 1 2 1 2 1 1 1 2 2 2 2 1 1 2 2 1 2 2 1 2 2 2 2 2 2
2 1 2 2 2 2 2 2 2 2 1 1 2 1}\\\hline
{\small 1 1 1 1 1 2 2 1 1 1 1 2 1 2 1 1 2 2 1 2 2 1 2 1 2 1 2 2 1 2 2 2 1 2 1
2 2 2 1 2 2 2 2 2 2 2 1 2 2}\\\hline
{\small 1 1 1 1 1 2 2 1 1 2 2 2 1 2 1 2 1 2 2 1 2 2 1 2 1 2 2 1 2 2 1 2 2 2 2
2 2 2 1 2 2 2 2 2 2 2 2 1 2}\\\hline
{\small 1 1 1 2 2 2 2 1 1 2 1 2 2 2 1 2 1 2 1 2 2 1 2 2 1 2 1 2 1 2 2 2 2 2 2
2 1 1 2 2 2 1 2 2 2 1 2 2 2}\\\hline
{\small 1 1 1 1 1 2 2 1 1 1 2 2 2 2 1 1 2 2 2 1 2 1 2 2 1 2 2 1 1 2 2 2 2 2 2
2 2 1 2 2 2 2 2 2 2 1 2 2 2}\\\hline
{\small 1 1 1 1 1 2 2 1 1 1 2 2 2 2 1 2 2 1 2 1 2 2 1 2 2 2 2 2 2 2 1 2 2 2 2
2 2 2 1 2 1 1 2 2 2 2 2 1 2}\\\hline
{\small 1 1 1 1 1 2 2 1 1 1 2 2 2 2 1 1 2 1 2 2 2 1 2 1 2 2 1 2 1 2 2 2 2 2 2
2 2 2 1 2 2 2 2 2 2 2 2 2 1}\\\hline
{\small 1 1 1 2 2 2 2 1 1 2 1 2 2 2 1 2 1 2 1 2 2 2 1 2 1 2 1 2 2 2 1 2 2 1 1
2 2 2 1 2 2 2 2 2 2 2 2 1 2}\\\hline
{\small 1 1 1 2 2 2 2 1 1 2 1 2 2 2 1 2 1 2 1 2 2 2 1 2 1 2 1 2 2 2 1 2 1 2 1
2 2 2 1 2 2 2 2 2 2 2 1 2 2}\\\hline
{\small 1 1 1 2 2 2 2 1 1 2 1 2 2 2 1 1 2 2 1 2 2 1 2 1 2 2 1 2 1 2 2 2 2 2 2
2 2 1 2 2 1 1 2 2 2 2 2 1 2}\\\hline
{\small 1 1 1 2 2 2 2 1 1 2 1 2 2 2 1 1 2 2 1 2 2 1 2 1 2 2 1 2 1 2 2 2 2 2 2
2 1 1 2 2 2 1 2 1 2 2 2 2 2}\\\hline
{\small 1 1 1 1 1 2 2 1 1 1 2 2 1 1 1 2 2 1 2 1 2 2 1 2 2 2 2 2 2 2 1 2 2 2 2
2 2 2 1 2 2 2 2 2 2 2 2 1 2}\\\hline
{\small 1 1 1 2 2 2 2 1 2 2 1 1 2 2 1 2 2 2 2 2 2 2 1 2 1 2 1 2 2 1 2 2 2 2 2
2 2 2 1 2 2 2 2 2 2 2 2 2 1}\\\hline
\end{tabular}
\
\]

\end{center}

\noindent No attempt has been made to gather matrices containing at least one
negative entry.

\section{Closing Words}

For any $3\times3$ unimodular matrix $M$ with $\alpha=2$, each entry of
$M^{-1}$ is a $2\times2$ minor of $M$, up to sign, hence has the form%
\[%
\begin{array}
[c]{ccc}%
\pm\left(  a\,d-b\,c\right)  , &  & \left\vert a\right\vert ,\left\vert
b\right\vert ,\left\vert c\right\vert ,\left\vert d\right\vert \leq2;
\end{array}
\]
thus%
\[
\left\vert a\,d-b\,c\right\vert \leq\left\vert a\,d\right\vert +\left\vert
b\,c\right\vert \leq4+4=8,
\]
i.e., $\beta\leq8$. \ More generally, for $n\times n$ unimodular $M$ with
$\alpha=2$, we have a theoretical bound $\beta\leq(n-1)!\,2^{n-1}$. \ Reason:
an $\left(  n-1\right)  \times\left(  n-1\right)  $ determinant with entries
in $\left[  -2,2\right]  $ is a sum of $\left(  n-1\right)  !$ signed
products, each of magnitude at most $2^{n-1}$. \ This is much larger than
empirical studies suggest is necessary. \ An additional requirement that
$M^{-1}$ contains no zeroes further reduces the bound on $\beta$, at least
initially:%
\[%
\begin{tabular}
[c]{|c|c|c|c|c|c|}\hline
$n$ & $3$ & $4$ & $5$ & $6$ & $7$\\\hline
$\beta_{\text{theor}}$ & $8$ & $48$ & $384$ & $3840$ & $46080$\\\hline
$\beta_{\text{empir}}$ & $6$ & $30$ & $182$ & $1122$ & -\\\hline
$\beta_{\text{zerofr}}$ & $5$ & $26$ & $182$ & $1122$ & $\geq7926$\\\hline
\end{tabular}
\ \ \ \ \ \ \
\]
The theoretical bound grows factorially, but empirical searches show that
maximal cofactors of unimodular matrices are substantially smaller. \ Imposing
the zerofree condition tightens the bound for small dimensions ($n=3,4$), but
for moderate dimensions ($n=5,6$) zerofree \& unrestricted empirical estimates
coincide, suggesting that the constraint ceases to be restrictive. \ [See the
Supplement.] \ The $7\times7$ scenario is beyond our computational reach; the
best matrix found so far is%
\[
\left(
\begin{array}
[c]{ccccccc}%
1 & 1 & 1 & 1 & 2 & 2 & 2\\
1 & 2 & -2 & -2 & 2 & 2 & 2\\
2 & 1 & -2 & -2 & -1 & -2 & -2\\
2 & 2 & 2 & 2 & -1 & 2 & -2\\
2 & 2 & 2 & -1 & 2 & -2 & 2\\
2 & -2 & 2 & -2 & -1 & 2 & -1\\
2 & -2 & -2 & 2 & 2 & -1 & 1
\end{array}
\right)
\]
and we expect that someone else can surely improve upon this.

\section{Addendum}

Let $\varphi$ be the Euler totient. \ Returning to $n=2$ and $\alpha=\beta$,
we obtain the following.\smallskip

\textit{Proposition 5}. \ Fix $k\geq2$. The number of matrices with $\alpha=k$
is $2\varphi(k)-1$.\smallskip

\textit{Proof}. \ Each such matrix $M$ can be represented as%
\[
\left(
\begin{array}
[c]{cc}%
a & b\\
c & d
\end{array}
\right)
\]
where $a,b,c,d$ are all positive and $d=k>\max\{a,b,c\}$. \ Let $\mathcal{M}$
denote the set of all $M$ and let $\varepsilon=a\,d-b\,c$; thus $\varepsilon
=\pm1$ and consequently $\gcd(b,k)=1$. \ Reason: if a positive integer $\ell$
divides both $b$ and $d$, then $\ell$ must divide $a\,d$ and $b\,c$, i.e.,
$\ell$ must divide $\varepsilon$, which implies that $\ell=1$. \ Therefore $b$
is invertible $\operatorname{mod}k$. (The same is also true for $c$.) \ Define
a function%
\[
\mathcal{M}\overset{f}{\longrightarrow}\{-1,+1\}\times\left(
\mathbb{Z}
/k%
\mathbb{Z}
\right)  ^{\ast}%
\]
by $f(M)=\left(  \varepsilon,b\right)  $. \ It suffices to demonstrate that
(i)\ $f$ is injective and (ii) exactly one codomain point is not in the image
of $f$. \ Starting with $\left(  \varepsilon,b\right)  \neq(-1,1)$, we have%
\[%
\begin{array}
[c]{ccc}%
c\equiv-\varepsilon\,b^{-1}\operatorname{mod}k, &  & a\,d\equiv\left(
\varepsilon+b\,c\right)  \operatorname{mod}k.
\end{array}
\]
For example, if $k=5$, then%
\[%
\begin{tabular}
[c]{|c|c|c|c|c|c|}\hline
$\left(  \varepsilon,b\right)  $ & $-\varepsilon\;$ & $\;b^{-1}\;$ &
$-\varepsilon\,b^{-1}$ & $\varepsilon+b\,c$ & $(a,c)$\\\hline
$(-1,4)$ & $1$ & $4$ & $4$ & $15$ & $(3,4)$\\\hline
$(1,3)$ & $-1$ & $2$ & $3$ & $10$ & $(2,3)$\\\hline
$(1,4)$ & $-1$ & $4$ & $1$ & $5$ & $(1,1)$\\\hline
$(-1,3)$ & $1$ & $2$ & $2$ & $5$ & $(1,2)$\\\hline
$(-1,2)$ & $1$ & $3$ & $3$ & $5$ & $(1,3)$\\\hline
$(1,2)$ & $-1$ & $3$ & $2$ & $5$ & $(1,2)$\\\hline
$(1,1)$ & $-1$ & $1$ & $4$ & $5$ & $(1,4)$\\\hline
\end{tabular}
\ \
\]
reproducing the seven-matrix result in Section 1. \ Starting with $\left(
\varepsilon,b\right)  =(-1,1)$, however, we would have $(a,c)=(0,1)$ but this
contradicts zerofreeness.

Returning to $n=3$, the sequence for $\alpha=3\leq\beta\leq15$ is%
\[
1,1,6,7,14,16,12,8,12,9,7,0,8
\]
and, to $n=4$, the sequence for $\alpha=3\leq\beta\leq105$ is%
\begin{align*}
&  163,183,380,393,771,853,1217,1182,1934,1720,2563,1826,2983,2606,3476,\\
&
1974,3471,3857,3559,2412,4114,2750,5309,2430,3590,3780,4030,2662,4136,3162,\\
&
3374,2105,4874,2879,3596,2128,2823,3348,3090,1850,2885,2563,3102,1431,2574,\\
&  1911,1999,1723,1878,1633,1907,983,2354,808,1484,754,1203,1842,1136,512,\\
&  690,606,1282,792,628,428,750,121,407,408,404,89,1041,112,80,418,148,84,\\
&  178,0,114,156,124,0,286,29,0,112,0,0,82,0,44,36,0,0,0,0,0,40,0,0,26.
\end{align*}
To do likewise for $n=5$ would necessitate examining $\alpha=3\leq\beta
\leq1023$. \ 

Returning to $(n,\alpha,\beta)=(6,2,3)$, we count $154$ matrices, of which six
are positive:%
\[%
\begin{array}
[c]{ccccc}%
\left(
\begin{array}
[c]{cccccc}%
1 & 1 & 1 & 1 & 1 & 2\\
1 & 1 & 1 & 1 & 2 & 1\\
1 & 1 & 1 & 2 & 1 & 1\\
1 & 1 & 2 & 1 & 2 & 2\\
1 & 2 & 1 & 2 & 1 & 2\\
2 & 1 & 1 & 2 & 2 & 1
\end{array}
\right)  , &  & i\left(
\begin{array}
[c]{cccccc}%
1 & 1 & 1 & 1 & 1 & 2\\
1 & 1 & 1 & 1 & 2 & 1\\
1 & 1 & 1 & 2 & 2 & 2\\
1 & 1 & 2 & 2 & 1 & 1\\
1 & 2 & 2 & 1 & 1 & 1\\
2 & 1 & 2 & 2 & 2 & 2
\end{array}
\right)  , &  & \left(
\begin{array}
[c]{cccccc}%
1 & 1 & 1 & 1 & 1 & 2\\
1 & 1 & 1 & 1 & 2 & 1\\
1 & 1 & 1 & 2 & 2 & 2\\
1 & 1 & 2 & 2 & 1 & 2\\
1 & 2 & 2 & 1 & 1 & 2\\
2 & 1 & 2 & 1 & 1 & 2
\end{array}
\right)  ,
\end{array}
\]%
\[%
\begin{array}
[c]{ccccc}%
\left(
\begin{array}
[c]{cccccc}%
1 & 1 & 1 & 1 & 1 & 2\\
1 & 1 & 1 & 1 & 2 & 1\\
1 & 1 & 1 & 2 & 2 & 2\\
1 & 1 & 2 & 2 & 1 & 2\\
1 & 2 & 2 & 1 & 1 & 2\\
2 & 1 & 2 & 2 & 2 & 2
\end{array}
\right)  , &  & \left(
\begin{array}
[c]{cccccc}%
1 & 1 & 1 & 1 & 1 & 2\\
1 & 1 & 1 & 2 & 2 & 2\\
1 & 1 & 2 & 1 & 2 & 2\\
1 & 2 & 1 & 2 & 2 & 1\\
1 & 2 & 2 & 2 & 2 & 2\\
2 & 2 & 2 & 1 & 2 & 2
\end{array}
\right)  , &  & \left(
\begin{array}
[c]{cccccc}%
1 & 1 & 1 & 2 & 2 & 2\\
1 & 2 & 2 & 1 & 1 & 2\\
1 & 2 & 2 & 2 & 2 & 2\\
2 & 1 & 2 & 1 & 2 & 1\\
2 & 1 & 2 & 2 & 2 & 2\\
2 & 2 & 2 & 1 & 2 & 2
\end{array}
\right)  .
\end{array}
\]
Results for the cases $(6,3,4)$ and $(6,3,5)$ would be good to see someday.

\section{Supplement}

Given $n\geq2$ and $\alpha\geq2$, define
\[
\gamma_{n}(\alpha)=\max\left\{  \beta\geq2:\exists\text{ }n\times n\text{
unimodular zerofree }M\text{ with }\left\Vert M\right\Vert =\alpha\text{,
}\left\Vert M^{-1}\right\Vert =\beta\right\}  .
\]
For example, $\gamma_{n}(2)=\beta_{\text{zerofr}}$ for $3\leq n\leq7$ as
described in Section 7. \ Because $\left\Vert M\right\Vert =\left\Vert
M^{-1}\right\Vert $ for $2\times2\ $matrices (Section 1), $\gamma_{2}(k)=k$
follows immediately.\smallskip

\textit{Conjecture 6}. \ We have%
\[
\gamma_{3}(k)=\left\{
\begin{array}
[c]{ccc}%
5 &  & \text{if }k=2,\\
k\left(  2k-1\right)  &  & \text{if }k\geq3;
\end{array}
\right.
\]%
\[
\gamma_{4}(k)=\left\{
\begin{array}
[c]{ccc}%
26 &  & \text{if }k=2,\\
k\left(  2k-1\right)  \left(  2k+1\right)  &  & \text{if }k\geq3;
\end{array}
\right.
\]%
\[%
\begin{array}
[c]{ccc}%
\gamma_{5}(k)=\left\{
\begin{array}
[c]{ccc}%
182 &  & \text{if }k=2,\\
1023 &  & \text{if }k=3,\\
3420 &  & \text{if }k=4;
\end{array}
\right.  &  & \gamma_{6}(2)=1122.
\end{array}
\]
Matrices
\[%
\begin{array}
[c]{ccc}%
\left(
\begin{array}
[c]{ccccc}%
1 & 1 & 2 & 3 & 3\\
2 & 4 & -5 & 5 & 5\\
2 & 5 & -5 & -5 & -5\\
3 & 5 & 5 & 5 & -4\\
4 & 5 & 4 & -5 & 5
\end{array}
\right)  , &  & \left(
\begin{array}
[c]{ccccc}%
1 & 1 & 1 & 2 & 3\\
1 & 6 & -6 & -5 & 6\\
1 & -6 & 6 & -6 & 5\\
4 & 6 & 5 & -6 & -6\\
5 & 6 & 6 & 6 & 6
\end{array}
\right)  ,
\end{array}
\]%
\[%
\begin{array}
[c]{ccc}%
\left(
\begin{array}
[c]{ccccc}%
1 & 6 & 7 & 7 & 7\\
3 & 7 & 7 & -7 & -7\\
3 & -5 & 2 & -2 & -7\\
5 & -7 & 7 & 7 & -6\\
7 & -7 & 7 & -6 & 6
\end{array}
\right)  , &  & \left(
\begin{array}
[c]{ccccc}%
2 & 4 & 5 & 5 & 5\\
7 & 7 & -6 & -8 & -8\\
7 & 8 & 7 & 7 & 8\\
8 & -7 & -2 & 8 & -8\\
8 & -7 & -7 & -8 & 8
\end{array}
\right)
\end{array}
\]
indicate that $\gamma_{5}(5)\geq8645$, $\gamma_{5}(6)\geq18282$, $\gamma
_{5}(7)\geq33033$, $\gamma_{5}(8)\geq54128$ and\pagebreak\
\[%
\begin{array}
[c]{ccc}%
\left(
\begin{array}
[c]{cccccc}%
1 & 1 & 1 & 1 & 2 & 2\\
1 & 2 & 3 & -3 & 3 & -3\\
1 & 3 & 3 & 3 & -3 & -2\\
2 & 3 & 3 & 2 & 3 & 3\\
2 & 3 & -3 & 3 & 3 & -3\\
3 & -3 & 3 & 3 & 2 & -3
\end{array}
\right)  , &  & \left(
\begin{array}
[c]{cccccc}%
1 & 1 & 1 & 1 & 3 & 4\\
1 & 4 & -4 & -4 & -4 & -2\\
2 & 3 & 4 & -4 & 4 & 4\\
2 & 4 & -4 & 4 & 3 & 4\\
3 & 4 & 4 & 3 & 4 & -4\\
4 & 3 & 4 & 4 & -3 & 3
\end{array}
\right)  ,
\end{array}
\]%
\[%
\begin{array}
[c]{ccc}%
\left(
\begin{array}
[c]{cccccc}%
1 & 1 & 2 & 3 & 4 & 5\\
2 & 4 & -5 & 4 & 5 & -5\\
3 & 5 & 5 & 5 & 4 & 4\\
4 & -5 & -4 & 1 & 5 & 4\\
5 & 3 & -5 & 2 & -5 & 5\\
5 & -5 & 5 & -5 & -4 & -5
\end{array}
\right)  , &  & \left(
\begin{array}
[c]{cccccc}%
1 & 2 & 3 & 3 & 5 & 5\\
4 & 4 & 6 & -6 & 3 & 6\\
5 & -6 & -4 & -6 & 5 & -3\\
6 & 5 & -4 & 6 & 6 & 2\\
6 & 6 & -6 & -5 & -6 & -5\\
6 & -4 & 3 & 6 & -5 & -4
\end{array}
\right)
\end{array}
\]
indicate that $\gamma_{6}(3)\geq9435$, $\gamma_{6}(4)\geq33964$, $\gamma
_{6}(5)\geq92069$, $\gamma_{6}(6)\geq201431$. \ Again someone else can surely
improve upon these.

It is surprising that, to the best of our knowledge, the conjectured formulas
for $3\leq n\leq4$ are new. \ To prove $\gamma_{3}(k)\geq k\left(
2k-1\right)  $ and $\gamma_{4}(k)\geq k\left(  2k-1\right)  \left(
2k+1\right)  $, let%
\[%
\begin{array}
[c]{ccc}%
U=\left(
\begin{array}
[c]{ccc}%
1 & 2 & 3\\
k & 1 & -k\\
k & k & k-1
\end{array}
\right)  , &  & V=\left(
\begin{array}
[c]{cccc}%
1 & k & k & k\\
k-2 & k & -k+1 & -k\\
k & 3 & 2 & -k+1\\
k & -k+1 & k & -k
\end{array}
\right)  ;
\end{array}
\]
these are unimodular zerofree matrices with inverses
\[
U^{-1}=\left(
\begin{array}
[c]{ccc}%
-k^{2}-k+1 & -k-2 & 2k+3\\
2k^{2}-k & 2k+1 & -4k\\
-k^{2}+k & -k & 2k-1
\end{array}
\right)  ,
\]%
\[
V^{-1}=\left(
\begin{array}
[c]{cccc}%
-8k^{2}+2k+1 & -2k^{3}-3k^{2}+k & 4k^{3}-k & -2k^{3}-k^{2}+2k\\
4k^{2}-7k & k^{3}-3k & -2k^{3}+3k^{2}+k & k^{3}-k^{2}-2k+1\\
4k^{2}+k+2 & k^{3}+2k^{2}+k+1 & -2k^{3}-k^{2}-k & k^{3}+k^{2}\\
-8k^{2}+14k-4 & -2k^{3}+5k-2 & 4k^{3}-6k^{2}+1 & -2k^{3}+2k^{2}+3k-3
\end{array}
\right)
\]
and the lower bounds follow. \ To prove equality, however, is an open challenge.

\section{Acknowledgements}

The creators of Mathematica earn my gratitude every day:\ this paper could not
have otherwise been written. \ For the first time, I\ have used the Microsoft
Copilot generative AI chatbot for assistance in writing/testing code. \ It
introduced me to the extraordinary accelerated performance associated with
CUDA programming on Nvidia GPUs.

\end{document}